\documentclass{amsart}
\usepackage{amssymb}
\usepackage{amsmath}
\usepackage{amsfonts}

\newtheorem{theorem}{Theorem}
\theoremstyle{plain}

\newtheorem{corollary}{Corollary}

\newtheorem{lemma}{Lemma}

\newtheorem{proposition}{Proposition}
\newtheorem{remark}{Remark}

\numberwithin{equation}{section}

\begin{document}
\title[Metric complexity]{A "metric" complexity for weakly chaotic systems}
\author{Stefano Galatolo}
\address{Dipartimento di Matematica Applicata, Universit\`{a} di Pisa, via
Buonarroti 1 Pisa}
\email{galatolo@dm.unipi.it}
\urladdr{http://www2.ing.unipi.it/$\sim $d80288}
\date{October 30, 2006}
\thanks{}

\begin{abstract}
We consider the number of Bowen sets which are necessary to cover a large
measure subset of the phase space. This introduce some complexity indicator
characterizing different kind of (weakly) chaotic dynamics. Since in many
systems its value is given by a sort of local entropy, this indicator is
quite simple to be calculated. We give some example of calculation in
nontrivial systems (interval exchanges, piecewise isometries e.g.) and a
formula similar to the Ruelle-Pesin one, relating the complexity indicator
to some initial condition sensitivity indicators playing the role of
positive Lyapunov exponents.
\end{abstract}

\maketitle

\section{Introduction}

Many techniques and results have been developed for the study of smooth
hyperbolic systems (systems where the dynamics is given by smooth functions
and distances between typical nearby initial conditions expand or contract
exponentially fast). In recent times, systems whose dynamics is not regular
(sometime discontinuous) and/or not hyperbolic (no exponential
contraction/expansion) are more and more important in various kind of
applications (interval exchanges, piecewise isometries: \cite{KWC}, \cite{AKT}
, Hamiltonian systems with stable islands: \cite{AG}, \cite{AZ}, symbolic
systems and automata such as substitutions and similar).

Such systems often have zero entropy, but their dynamics is far to be simple
and predictable because there is still a "weak" initial condition
sensitivity "slowly" separating nearby starting orbits.

The need to provide complexity indicators which can describe and quantify
this "weakly"\ chaotic behavior lead in the mathematical literature to many
definitions and different notions of \emph{complexity} (or generalized
entropies).

The first natural attempt is to repeat the same construction leading to the
K-S entropy (considering first a partition of the space, considering the
induced symbolic system, and so on...) replacing the usual formula for the
Shannon entropy of a symbolic system ($-\sum p_{i}\log p_{i}$ ) with a
different one (here the physical literature is huge, but there are few
rigorous results, see e.g. \cite{Bl2},\cite{Ta},\cite{Ts}). This kind of
construction often has the problem that the resulting indicator is not
continuous with respect to change of partitions (see, \cite{Bl2},\cite{Ta})
thus its physical meaning is compromised and the calculation of the suprema
over all partition is difficult.

To overcome this difficulty a more refined construction can be performed (%
\cite{feren}, \cite{KT}). This lead to a more stable definition and to an
invariant which can be calculated and has nontrivial values on interesting
examples. This indicator works in the measure-theoretic, ergodic framework
and is invariant under measure preserving transformations.

Another, topological approach considers the number of essentially
different orbits (orbits whose distance at a certain time is greater
than a given resolution $\epsilon )$ which appear in the system
(\cite{BHM},\cite{Ga3}, see also \cite{AG}, \cite{AZ} for many
variants on this theme) and consider how this number increases with
time. This lead to topological complexities which generalizes the
topological entropy. The disadvantage of a purely topological
approach can be understood comparing topological entropy with
Kolmogorov-Sinai entropy. The presence of a physical invariant
measure in system gives more weight to the most frequent (and
physically relevant) configurations, neglecting the least relevant
ones, which on the other hand are not neglected in the topological
approach.

Another approach to define complexity is to consider the complexity of
single orbits of the system (see e.g. \cite{Br} \cite{Ga3}), this complexity
indicator is then local, the global behavior can be given by the complexity
of a typical orbit, or the average with respect to some invariant measure.
The orbit complexity is given by the amount of information (algorithmic
information) which is necessary to describe the orbit up to some give
accuracy. If the accuracy is given by some partition or by an open cover the
notion is more measure theoretic or topology oriented. In this approach the
complexity indicators can be easily calculated in many interesting examples,
and there are connections with many other features of chaotic dynamics, such
as dimension of attractors and so on (see, for example \cite{Z}, \cite{Ga1},
\cite{BG}, \cite{GW}).

In this paper we follow an approach which defines a global indicator of
complexity and which is not only topological or measure theoretic. We will
define some (more rigid) indicators which are invariant under morphisms
which are both continuous and measure preserving.

Many interesting physical coordinate change are continuous and they preserve
some physical measure (for example if we are observing and reconstructing a
system trough some continuous observable, as in the nonlinear time series
framework, see e.g. \cite{KaS}, \cite{OY}).

We will construct a complexity indicator which is invariant for this kind of
morphisms and it is easy to be calculated. Moreover it has connections with
the other features of chaos.

Roughly speaking we will consider the number of "important", essentially
different orbits\ which appear in the system. The importance will be given
by the measure $\mu .$ More precisely, we will consider the number of Bowen
sets which are necessary to cover a large part of $\mu $ \ and we will
consider how this number increases with time. We will see that under mild
assumptions, this indicator is equivalent to the rate of decreasing of the
measure of a typical Bowen set (a sort of extension of the Brin-Katok theorem \cite{BK}).
This will allow an easy calculation of the complexity indicator in
nontrivial cases, as interval exchanges, piecewise isometries, the logistic
map, and some more examples, which are listed in section 3. In section 4 we
will consider a set of numbers \ describing the geometrical features of the
Bowen set, these numbers plays the role of the Lyapunov exponents,
describing initial condition sensitivities at different directions and
allowing a result similar to the Ruelle-Pesin formula.

\section{A "metric" complexity}

We consider a system $(X,T,\mu )$ of the following type: $X$ is a metric
space equipped with a distance $d$. The dynamics is given by a Borel map  $T:X\rightarrow X$  and $\mu $ is invariant for $T$.

Let us consider the Bowen set
\begin{equation*}
B(n,x,\epsilon )=\{y\in X:d(T^{i}(y),T^{i}(x))\leq \epsilon \ \forall i\
s.t.\ \ 0\leq i\leq n\}.
\end{equation*}
$B(n,x,\epsilon )$ is the set of points \textquotedblleft
following\textquotedblright\ the orbit of $x$ for $n$ steps at a distance
less than $\epsilon $. As the nearby starting orbits of $(X,T)$ diverges the
set $B(n,x,\epsilon )$ will be smaller and smaller as $n$ increases. If two
points are in the same set we can think that their orbits are similar (up to
a resolution given by $\epsilon ,$ for $n$ steps) if two points are in
different sets, their orbits are essentially different\footnote{%
Counting the number of essentially different orbits needed to cover the
whole space $X$ leads to the notion of \ topological entropy and to its
generalizations which can be called topological complexity of a system.}.

We want to consider the number of Bowen sets which is necessary to cover a
large (according to the measure $\mu $ ) part of $X$. This counts how many
different "important" orbits appears in the system. Here the notion of
importance if provided by the measure $\mu ,$ which will give different
weight to different parts of $X.$ This complexity depends both on the
metric, and ergodic features of the system and the notion is physically
relevant when we consider a physical invariant measure. Hence this notion is
related to the metric of the system (which induces the Lesbegue measure,
which induces the physical measure, see e.g. \cite{Y}) for this reason we
call it "metric complexity".

Let us hence consider the following
\begin{equation}
N(n,\epsilon ,\epsilon ^{\prime })=\min (\{k\in \mathbb{N}|\exists
x_{1},...,x_{k},\mu (\cup _{0\leq i\leq k}B(n,x_{i},\epsilon ))\geq
1-\epsilon ^{\prime }\})\label{N}
\end{equation}%
that is the number of Bowen sets that is necessary to cover a subset of $X$
whose measure is bigger than $1-\epsilon ^{\prime }.$ We want to consider the
asymptotic growing rate of $N(n,\epsilon ,\epsilon ^{\prime })$ as $n$
increases, when $\epsilon $ and $\epsilon ^{\prime }$ are small.

To formalize this, for each monotonic function $f(n)$ with $%
\mathrel{\mathop{lim}\limits_{n\rightarrow \infty }}f(n)=\infty $ we define
an indicator by comparing the asymptotic behavior of $\log (N(n,\epsilon
,\epsilon ^{\prime }))$ with $f$ \footnote{%
From now on, in the definition of indicators $f$ is always assumed to be
monotonic and tends to infinity.}. Hence let us consider%
\begin{equation*}
h_{\epsilon ,\epsilon ^{\prime }}^{f}(X,T,\mu )=\underset{n\rightarrow
\infty }{\lim \sup }\frac{\log (N(n,\epsilon ,\epsilon ^{\prime }))}{f(n)}
\end{equation*}%
this quantity is monotonic in $\epsilon $ and $\epsilon ^{\prime }$ and
hence we can consider the limits

\begin{equation*}
h^{f}(X,T,\mu )=\underset{\epsilon ^{\prime }\rightarrow 0}{\lim }\underset{%
\epsilon \rightarrow 0}{\lim }h_{\epsilon ,\epsilon ^{\prime }}^{f}(X,T,\mu
).
\end{equation*}
We will see (see proposition \ref{bkg}) that when $f$ is the identity ( $%
f(n)=n$), the quantity $h^{id}(X,T,\mu )$ equals the Kolmogorov-Sinai
entropy for a large family of systems.

Let us now consider the invariance properties of $h^{f}$ under isomorphisms
of systems.

\begin{theorem}
If $(X,T,\mu )$, $(Y,T^{\prime },\mu ^{\prime })$ are dynamical systems over
compact metric spaces $(X,d),(Y,d^{\prime }).$ Let $\phi $ be a measure
preserving hoemorphism such that the following diagram%
\begin{equation*}
\begin{array}{rcccl}
\  & \  & \phi & \  & \  \\
\  & X & \rightarrow & Y\  &  \\
T & \downarrow & \  & \downarrow & T^{\prime } \\
\  & X & \rightarrow & Y\  &  \\
\  & \  & \phi & \  & \
\end{array}%
\end{equation*}%
commutes, then $h^{f}(X,T,\mu )=h^{f}(Y,T^{\prime },\mu ^{\prime }).$
\end{theorem}

\emph{Proof.} Let us call $N_{1}(n,\epsilon ,\epsilon ^{\prime })$ the
number of Bowen sets that is necessary to cover a large subset of $X$ as
above, and $N_{2}(n,\epsilon ,\epsilon ^{\prime })$ be the number of Bowen
sets that is necessary to cover a large subset of $Y$. Since the spaces are
compact and $\phi $ is continuous then it is uniformly continuous. Let $g:%
\mathbb{R\rightarrow R}$ \ such that $d(x_{1},x_{2})\leq g(\epsilon )$ (with
$x_{i}\in X$) implies $d^{\prime }(\phi (x_{1}),\phi (x_{2}))\leq \epsilon .$
For each $n$ it holds $f(B(x,n,g(\epsilon )))\subset B(f(x),n,\epsilon )$
then let us suppose that $\{B(n,x_{1},g(\epsilon )),...,B(n,x_{k},g(\epsilon
))\}$ is a minimal cover of a large set $A\subset \cup _{0\leq i\leq
k}B(n,x_{i},g(\epsilon ))$ with measure $\mu (A)=1-\epsilon ^{\prime }$,
this implies that $f(A)\subset \cup _{0\leq i\leq k}B(n,f(x_{i}),\epsilon ).$
We recall that $\mu (A)=\mu ^{\prime }(f(A)).$ Hence $N_{2}(n,\epsilon
,\epsilon ^{\prime })\leq N_{1}(n,g(\epsilon ),\epsilon ^{\prime }).$ This
implies that $h_{g(\epsilon ),\epsilon ^{\prime }}^{f}(X,T,\mu )\geq
h_{\epsilon ,\epsilon ^{\prime }}^{f}(Y,T^{\prime },\mu ^{\prime })$ and $%
h^{f}(X,T,\mu )\geq h^{f}(Y,T^{\prime },\mu ^{\prime }).$ Similarly we can
prove the reverse inequality.$\square $

It is useful to consider a version of the Brin-Katok local entropy (\cite{BK}%
): let us define%
\begin{equation*}
\overline{BK}^{f}(x,\epsilon )=\mathrel{\mathop{limsup}\limits_{n\rightarrow
\infty}}\frac{-log(\mu (B(n,x,\epsilon )))}{f(n)},\underline{BK}%
^{f}(x,\epsilon )=\mathrel{\mathop{liminf}\limits_{n\rightarrow \infty}}%
\frac{-log(\mu (B(n,x,\epsilon )))}{f(n)}
\end{equation*}%
\begin{equation*}
\overline{BK}^{f}(x)=\mathrel{\mathop{lim}\limits_{\epsilon\rightarrow 0}}%
\overline{BK}^{f}(x,\epsilon ),\underline{BK}^{f}(x)=\mathrel{\mathop{lim}%
\limits_{\epsilon\rightarrow 0}}\underline{BK}^{f}(x,\epsilon ).
\end{equation*}%
When $f(n)=n$ is the identity then $BK^{f}$ is the Brin-Katok local entropy.
In \cite{BK} it is proved that if the system is ergodic $\overline{BK}%
^{id}(x)=\underline{BK}^{id}(x)=h_{\mu }(T)$ (the K-S entropy) for almost
each $x\in X$. Hence $\overline{BK}^{id}(x)$ and $\underline{BK}^{id}(x)$
are almost everywhere equal and they are invariant under $T$ .

In the general case however the invariance under $T$ holds under some mild
conditions

\begin{proposition} \label{prop}If $T$ is such that
\begin{itemize}
\item i) Almost each point $x$ has a small neighborhood $U$ such that $%
T|_U:U\rightarrow T(U)$ is an homeomorphism

\item ii) For each measurable $A$ it holds $\mu (T(A))\leq K\mu (A)$ for
some fixed constant $K$
\end{itemize}

\noindent then $\overline{BK}^{f}(x)=\overline{BK}^{f}(T(x))$ and $%
\underline{BK}^{f}(x)=\underline{BK}^{f}(T(x))$ for $\mu $ almost each $x$.
\end{proposition}

\emph{Proof.} First let us notice that
\begin{equation*}
B(n,x,\epsilon )=B(x,\epsilon )\cap T^{-1}(B(n-1,T(x),\epsilon ))
\end{equation*}%
then it is clear that (T preserves $\mu $) $\mu ({B(n,x,\epsilon ))}\leq \mu
(B(n-1,T(x),\epsilon ))$ and then $\overline{BK}^{f}(x)\geq \overline{BK}%
^{f}(T(x))$ and $\underline{BK}^{f}(x)\geq \underline{BK}^{f}(T(x))$ .

For the other inequality, we have that $T$ is a.e. a local homeomorphism,
let $x$ be a typical point and $\epsilon ^{\prime }<\epsilon $ such that $%
B(T(x),\epsilon ^{\prime })\subset T(B(x,\epsilon ))$. Obviously $%
B(n-1,T(x),\epsilon ^{\prime })\subset B(T(x),\epsilon ^{\prime })$. Now $%
B(n-1,T(x),\epsilon ^{\prime })\subset T(B(n,x,\epsilon ))$, this is true
because if $y\in B(n-1,T(x),\epsilon ^{\prime })$ then there is some $z\in
B(x,\epsilon )$ with $T(z)=y$. Now, if $z$ is such that $d(x,z)<\epsilon $, $%
T(z)\in B(n-1,T(x),\epsilon ^{\prime })$ with $\epsilon ^{\prime }<\epsilon $
then $d(T^{i}(z),T^{i}(x))<\epsilon $ for each $0\leq i\leq n$ and then $%
z\in B(n,x,\epsilon )$. By ii) $\frac{1}{K}\mu (B(n,x,\epsilon ))\geq \mu
(T(B(n,x,\epsilon )))\geq \mu (B(n-1,T(x),\epsilon ^{\prime }))$, and then $%
BK(x)\leq BK(T(x))$. $\Box $

The relation between $BK^{f}$ and $h^{f}$ in general is quite natural

\begin{proposition}
\label{bkg} If $\overline{BK}^{f}(x)=\underline{BK}^{f}(x)=BK^{f}$ almost
everywhere then
\begin{equation*}
BK^{f}=h^{f}(X,T,\mu ).
\end{equation*}
\end{proposition}

\begin{proof}
Since $\overline{BK}(x)=\underline{BK}(x)$ almost everywhere, for each $%
\varepsilon >0$ there is an $\overline{n}$ \ and a set $A_{\varepsilon ,%
\overline{n}}$ such that for each $n\geq \overline{n}$ and $x\in
A_{\varepsilon ,\overline{n}}$%
\begin{equation*}
\mu (B(n,x,\epsilon ))\leq 2^{-(BK(x)g(\epsilon )-\varepsilon )f(n)}
\end{equation*}
for some $g,$ such that $g(\epsilon )\rightarrow 1$ as $\epsilon \rightarrow
0$. Moreover the sequence $A_{\varepsilon ,\overline{n}}$ is increasing as $%
\overline{n}$ increases and $\mu (A_{\varepsilon ,\overline{n}})\rightarrow
1 $ as $\bar{n}\rightarrow \infty $. Let us fix an arbitrary small $%
\varepsilon $ and $\overline{n}$ such that
\begin{equation*}
\mu (A_{\varepsilon ,\overline{n}})\geq \frac{3}{4}.
\end{equation*}%
Now, let us consider $n\geq \overline{n}$ and a set $\{B(n,x_{1},\frac{%
\epsilon }{2}),...,B(n,x_{k},\frac{\epsilon }{2})\}$ covering a big subset
of $X$ as in the definition of $h^{f}(X,T,\mu )$. More precisely, we can
suppose that $\mu (\cup _{0\leq i\leq k}B(n,x_{i},\frac{\epsilon }{2}))\geq
\frac{3}{4}$ and hence $\mu (\cup _{0\leq i\leq k}B(n,x_{i},\frac{\epsilon }{%
2}))\cap A_{\varepsilon ,\overline{n}}\geq \frac{1}{2}.$

Now we remark that if $B(n,x_{i},\frac{\epsilon }{2})\cap A_{\varepsilon ,%
\overline{n}}\neq \emptyset $ then there is $x\in A_{\varepsilon ,\overline{n%
}}$ such that $B(n,x_{i},\frac{\epsilon }{2})\subset B(n,x,\epsilon ),$
hence $\mu (B(n,x_{i},\frac{\epsilon }{2}))\leq \mu (B(n,x,\epsilon ))\leq $
$2^{-(BK(x)g(\epsilon )-\varepsilon )f(n)}.$ Since each one of these sets $%
B(x,n,\frac{\epsilon }{2})$ has small measure and their union has measure
greater than $\frac{1}{2}$ then its number must be greater than $%
2^{(BK(x)g(\epsilon )+\varepsilon )f(n)-1}$ giving that $h_{\frac{\epsilon }{%
2},\frac{3}{4}}^{f}(X,T,\mu )\geq BK^{f}(x,\epsilon )$ almost everywhere,
hence $h^{f}(X,T,\mu )\geq BK^{f}(x,\epsilon )$ a.e.

For the other inequality, similar as before for each $\varepsilon $ there is
an $\overline{n}$ and a set $B_{\varepsilon ,\overline{n}}$ such that for
each $n\geq \overline{n}$ and $x\in B_{\varepsilon ,\overline{n}}$ it holds $%
\mu (B(n,x,\epsilon ))\geq 2^{-(BK(x)g(\epsilon )+\varepsilon )f(n)}$ for
some $g,$ such that $g(\epsilon )\rightarrow 1$ as $\epsilon \rightarrow 0$
and $\mu (B_{\varepsilon ,\overline{n}})\rightarrow 1$ as $\bar{n}%
\rightarrow \infty $. Let us consider $C=\{B(n,x_{1},\epsilon
),...,B(n,x_{k},\epsilon )\}$ such that $C$ is made of disjoint Bowen sets,
each $x_{i}$ is contained in $B_{\varepsilon ,\overline{n}}$ and $C$ is
maximal, in the sense that $\forall x\in B_{\varepsilon ,\overline{n}}$ then
$B(n,x,\epsilon )\cap B(n,x_{i},\epsilon )\neq \emptyset $ for some $%
B(n,x_{i},\epsilon )\in C.$ The set $C$ is finite because by definition of $%
B_{\varepsilon ,\overline{n}}$ each $B(n,x_{i},\epsilon )\in C$ has a
measure greater than $2^{-(BK(x)g(\epsilon )-\varepsilon )f(n)}$ and their
total measure must be less than $1.$ Thus the number of such set is less or
equal than $2^{(BK(x)g(\epsilon )-\varepsilon )f(n)}.$ Now we remark that if
$C$ is \ as before, then $C_{2}=\{B(n,x_{1},2\epsilon
),...,B(n,x_{k},2\epsilon )\}$ is a cover of $B_{\varepsilon ,\overline{n}}.$
Then we proved that there is a cover of some big as wanted subset (with
measure let us say, greater than $1-\epsilon ^{\prime }$ ) of $X$ made with
no more than $2^{(BK(x)g(\epsilon )+\varepsilon )f(n)}$ Bowen sets and this
proves $h_{\epsilon ^{\prime },2\epsilon }^{f}(X,T,\mu )\leq
BK^{f}(x,\epsilon )$ and $h^{f}(X,T,\mu )\leq BK^{f}(x)$ a.e.
\end{proof}

\begin{remark}
\label{rem2}If $\overline{BK}^{f}\geq \overline{BK}^{f}(x)\geq \underline{BK}%
^{f}(x)\geq \underline{BK}^{f}$ almost everywhere\footnote{%
This happen for example if $(X,T,\mu )$ is ergodic and it satisfies the
assumptions of proposition \ref{prop}.}, the above proof gives that%
\begin{equation*}
\underline{BK}^{f}\leq h^{f}(X,T,\mu )\leq \overline{BK}^{f}.
\end{equation*}
\end{remark}
For a natural example where $\overline{BK}^{f}(x)>\underline{BK}^{f}(x)$
a.e. see section \ref{exa}.

Proposition \ref{bkg} allows to easily calculate $h^{f}(X,T,\mu ).$ If the
assumptions of the proposition are verified, instead to construct a global
cover of the system by Bowen sets we only need to look the behavior of the
measure of a typical Bowen set. To give an example of nontrivial
calculation, in next section we calculate the complexity of typical Interval
Exchange Transformations, the Logistic map the Feigenbaum point, the
Casati-Prosen map.

\section{Some example}

As said before, since the assumptions of Proposition \ref{bkg} are
mild and easy to be verified we can apply it in many cases and
estimate $ h^{f}(X,T,\mu )$ by $BK^{f}(x,\epsilon ),$ which is an
estimation of initial condition sensitivity at typical points. We
give some example of this application on some non trivial examples.

\subsection{General piecewise isometries}

Let consider a nontrivial family of systems for which we can have an
upper estimation for the complexity. Piecewise Isometries (PI)\ are
simple families of dynamical systems that show dynamical complexity
while not being hyperbolic in any senses; classical examples in one
dimension are, interval exchange transformations (IETs, see also
below). PIs have also been found to arise in several applications
such as in digital filter models and billiard systems ( see
\cite{ash}, \cite{Go}).

It is conjectured that the symbolic dynamics of a PI has polynomial
complexity (in the sense that the number of different names of
subcilynders appearing in the dynamics grow polinomially with the
length, for some works on this direction see e.g. \cite{AKT},
\cite{Bu}, \cite{Ka}). We give an upper bound of our definition of
complexity. This correspond to a polinomial bound on the growth of
Bowen sets necessary to cover the invariant measure (instead of
cylinders).

Let us recall briefly the class of systems we are considering. Let $X=%
\mathbb{R}^{n},$ Let us suppose that $P_{1},...,P_{m}$ is a measurable
partition of $X.$

A piecewise isometry $T:X\rightarrow X$ is a map defined in the following
way: let $A_{1},...,A_{m}:X\rightarrow X$\ be a set of isometries, then $%
T(x)=A_{i}(x)\Longleftrightarrow x\in P_{i}$. The sets $P_{i}$ are called
atoms and most of the literature consider piecewise linear atoms. We will
consider a more general situation.

In our piecewise isometries, the only source of initial condition
sensitivity is the presence of discontinuities at the boundary of atoms. Let
$Y=\cup _{i\leq m}\partial P_{i}$. If for each $i\leq n$ it holds $%
d(T^{i}(x),Y)\geq r$ then we know that the Bowen set satisfies $%
B(x,n,\epsilon )\supset B_{r}(x)${}for each $\epsilon >r.$ Hence the initial
condition sensitivity depends on the speed a typical orbit approaches the
discontinuity set $Y$.

To estimate this we will use the following simple result (\cite{DeGa} Lemma
2): given $Y\subset X$, let us define the $r$ neighborhood of $Y$ by
\begin{equation*}
B_{r}(Y)=\{x\in X,\,d(x,Y)<r\}
\end{equation*}%
and consider \underline{$d$}$_{\mu }(Y)=%
\mathrel{\mathop{{\liminf}}\limits_{\epsilon \rightarrow
0}}\,\,\frac{\log (\mu (B_{\epsilon }(Y)))}{\log (\epsilon )}.$ We remark
that if $Y=x$ is a point this gives the definition of lower local dimension
of $\mu $ at $x.$We recall that if \underline{$d$}$_{\mu }(x)=%
\mathrel{\mathop{{\limsup}}\limits_{\epsilon \rightarrow
0}}\,\,\frac{\log (\mu (B_{\epsilon }(x)))}{\log (\epsilon )}=d_{\mu }(x)$
this is called the local dimension of $\mu $ at $x.$

\begin{lemma}
\label{GAN2} Let $(X,T,\mu )$ be a measure preserving transformation, $%
Y\subset X$. If $\alpha >\frac{1}{\underline{d}_{\mu }(Y)}$ then for almost
each $x\in X$:
\begin{equation*}
\mathrel{\mathop{\liminf}\limits_{n\rightarrow \infty}}n^{\alpha
}d(T^{n}(x),Y)=\infty .
\end{equation*}
\end{lemma}

Hence we obtain the following

\begin{proposition}
If $T$ is an ergodic piecewise isometry as defined above, $Y=\cup
_{i\leq m}\partial P_{i}$ and $\underline{d}_{\mu }(Y)\neq 0$,
moreover if the local
dimension $d_{\mu }(x)$ is well defined and a.e. constant on $X,$ then%
\begin{equation*}
h_{\mu }^{\log }(T)\leq \frac{d_{\mu }(x)}{\underline{d}_{\mu }(Y)}
\end{equation*}
\end{proposition}

\begin{proof}
First we remark that since $d=\overline{d}_{\mu }(Y)\neq 0$ then
$\mu (Y)=0.$ This, together with the other properties of piecewise
isometries implies that $T$ satisfies the assumptions of proposition
\ref{prop}, hence by remark \ref{rem2} it is sufficient to estimate
the behavior of $\mu (B(x,n,\epsilon )).$ First we remark that we can suppose $%
\mathrel{\mathop{\liminf}\limits_{n\rightarrow
\infty}}d(T^{n}(x),Y)=0,$ otherwise the statement is trivial (because the
typical orbit never approaches to the discontinuity). In this case, as
remarked above, by Lemma \ref{GAN2} we have that for almost each $x\in
X$, small $\varepsilon >0$ it holds $B(x,n,\epsilon )\supset B_{n^{\frac{-1}{%
d+\varepsilon }}}(x)$ eventually with respect to $n.$ Then if $n$ is big
enough $\mu (B(x,n,\epsilon ))\geq \mu (B_{n^{\frac{-1}{d+\varepsilon }%
}}(x)).$ By the assumptions on the local dimension of the system
then we have that again, if $n$ is big enough, if $\varepsilon
,\epsilon ^{\prime }$ are small $\mu (B(x,n,\epsilon ))\geq
n^{(d_{\mu }(x)-\epsilon ^{\prime })(\frac{1}{d+\varepsilon })}.$
Which gives the statement.
\end{proof}

\subsection{Interval Exchanges}

Interval Exchanges are close relatives of surface flows, these maps are
particular bijective piecewise isometries of the unit interval, whose atoms
are intervals and which preserve the Lesbegue measure. In this section we
apply a \ result of Boshernitzan about a full measure class of uniquely
ergodic interval exchanges to estimate their metric complexity. We refer to
\cite{Bo2} for generalities on this important class of maps.

Let $T$ be some interval exchange. Let $\delta (n)$ be the minimum distance
between the discontinuity points of $T^{n}.$ We say that $T$ has the
property $\tilde{P}$ if there is a\ constant $C$ and a sequence $n_{k}$ such
that $\delta (n_{k})\geq \frac{C}{n_{k}}.$

\begin{lemma}
\label{bosh9}(by \cite{Bo2}) The set of interval exchanges having the
property $\tilde{P}$ has full measure in the space of interval exchange maps$%
.$
\end{lemma}

From Lemma \ref{GAN2} it easily follows that

\begin{corollary}
\label{dist}For each interval exchange $T$ and each $\epsilon >0$, for
almost each $x$ the distance from the orbit of $x$ to the discontinuity set
of $T$ is estimated as follows. If $y_{1},...,y_{k}$ are the discontinuity
points of $T$ then eventually with respect to $n$
\begin{equation*}
\underset{i\leq n,j\leq k}{\min }d(T^{i}(x),y_{j})>n^{-1-\epsilon }.
\end{equation*}
\end{corollary}

Since the initial condition sensitivity of interval exchanges is
determined by the speed of approaching of starting points to the
discontinuities, these results will allow to estimate $h_{\mu
}^{f}(T).$ Indeed by the above corollary \ref{dist} we know that if
$T$ is ergodic, for almost each $x$,
for each $\varepsilon >0$ eventually $\mu (B(n,x,\varepsilon ))\geq $ $%
n^{-1-\epsilon }.$ Since an interval exchange satisfies the assumptions of
remark \ref{rem2} then this implies that $h_{\mu }^{\log }(T)\leq 1.$

On the other hand the converse estimation follows from the remark that if $%
x_{0}$ is a discontinuity point, and%
\begin{equation*}
\underset{i\leq n,T^{i}(x)\leq x_{o}}{\min }d(T^{i}(x),x_{0})=l_{1}(n)\ and%
\underset{i\leq n,T^{i}(x)\geq x_{o}}{\min }d(T^{i}(x),x_{0})=l_{2}(n)
\end{equation*}%
(the minimum distance after $n$ steps of the orbit on the left and on the
right side of the discontinuity $x_{0}$) then for small $\epsilon ,$ $%
B(n,x,\epsilon )\subseteq (x-l_{1}(n),x+l_{2}(n)).$ Now we have to estimate
from above the speed of approaching to the discontinuity on both sides.
Using property $\tilde{P},$ like in \cite{Ga1} we can obtain the following

\begin{proposition}
Let $T$ be an IET with property $\tilde{P}$ as before, then $h_{\mu }^{\log
}(T)\geq 1.$

\begin{proof}
If $T$ has $m$ discontinuity points, $T^{n}$ has $nm$ discontinuity points
and they will divide the unit segment into $nm+1$ small segments. The total
length is $1$, then among these small segments there are at least $\frac{nm}{%
2}$ ones with length less or equal than $\frac{2}{mn+1}.$ Let us denote by $%
J_{n}$ the union of these segments. By property $\tilde{P}$ there is a
sequence $n_{k}$ such that the segments in $J_{n_{k}}$ are longer than $%
\frac{C}{n_{k}},$ by this $\mu (J_{n_{k}})\geq \frac{mC}{2}$. Hence there is
a set $B$ with positive measure, $\mu (B)\geq \frac{mC}{2}$ such that if $%
x\in B$ then $x$ is contained in infinitely many $J_{n_{k}}.$ Let us
notice at this point that if $x\in J_{n_{k}}$ then the
discontinuities of $T^{n_{k}}$ near $x$ are the ends of the small
interval $(y_{i},y_{j})\subset
J_{n_{k}}$ containing $x,$ hence for small $\epsilon $ the Bowen set around $%
x$ satisfies $B(n_{k}+1,x,\epsilon )\subseteq (y_{i},y_{j}).$ Recalling that
$\mu (J_{n_{k}})\geq \frac{mC}{2}$ now, we estimate (see eq. \ref{N}) $N(n_{k}+1,\epsilon ,1-%
\frac{mC}{4}).$ To cover a set with measure greater than
$1-\frac{mC}{4}$ we need to cover at least half of $J_{n_{k}}$, but
his intervals (and respective Bowen sets) have measure less or equal
than $\frac{2}{mn_{k}+1},$ hence we need at least
$\frac{mn_{k}+1}{2}\frac{mC}{4}$ sets, which gives the statement.
\end{proof}
\end{proposition}

Collecting the above results we have the following estimation of the
complexity for typical interval exchanges.

\begin{proposition}
If $T$ is an IET with property $\tilde{P}$ then $h_{\mu }^{\log }(T)=1.$
\end{proposition}

The situation for nontypical IET in general is much more
complicated. We expect arithmetical phenomena like in section
\ref{exa} to happen.

\subsection{Casati Prosen map}

In this subsection we will consider the Casati Prosen map, the map acts on
the unit square, is weakly chaotic and it is not a piecewise isometry. This
kind of map was introduced by Casati and Prosen \cite{CP2} in connection
with the mixing properties of flows in certain triangular billiards \cite%
{CP1}. We will give an upper estimation of its complexity.

Let us define the map: let $\theta (q)$ be the discontinuous function over
the circle given by $\theta (q)=-1$ if $0\leq q\leq 1/2$ and $\theta (q)=1$
otherwise.

For any $\alpha ,\beta \in \lbrack 0,1]$, we define the map $T_{\alpha
,\beta }$ as
\begin{equation*}
T_{\alpha ,\beta }(q,p)\,=\,(q+p+\beta \,,\,p+\alpha \,\theta (q))\,\,\,\,%
\hbox{\rm mod\ }1.
\end{equation*}%
We remark that $T_{\alpha ,\beta }$ can be written as the composition of
three elementary maps,
\begin{equation*}
T_{\alpha ,\beta }\,=\,B\circ R\circ G_{\alpha },
\end{equation*}%
where $B$ is represented by the matrix $\left(
\begin{array}{cc}
1 & 1 \\
0 & 1%
\end{array}%
\right) $ (a skew translation), $R(q,p)=(q+\beta ,p)$ is a translation in
the $q$ direction and $G_{\alpha }$ is the \textit{discontinuous part} of
the dynamics $G_{\alpha }(q,p)=(q,p+\alpha \,\theta (q))$ this discontinuous
map cuts the square along the lines $\rho \cup \rho ^{\prime }=\left(
\{1/2\}\times \lbrack 0,1[\right) \cup \left( \{0\}\times \lbrack
0,1[\right) $. translating the two pieces in opposite directions along the
line. Hence initially close orbits separate in a way that the distance
increases linearly with time by the skew translation, until they are
drastically separated by the discontinuity. The Lesbegue measure $\lambda $
is invariant for the map. It is surprising that there are few rigorous
results about ergodic properties of such map. As far as we know, even
ergodicity for $\alpha \neq 0$ is still not proven (even if probably true
for irrational values for the parameters). The map satisfies the assumptions
of remark \ref{rem2}, hence to give an estimation of the complexity it is
sufficient to estimate the behavior of $\mu (B(x,n,\epsilon )).$

\begin{proposition}
If $(X,T_{\alpha ,\beta },\lambda )$ is the Casati Prosen map then $h^{\log
}(T_{\alpha ,\beta })\leq 3.$
\end{proposition}

\begin{proof}
Let us consider $Y=\rho \cup \rho ^{\prime }$ (the discontinuity set) since
we consider the Lesbegue measure we have $\underline{d}_{\lambda }(Y)=1$,
hence by lemma \ref{GAN2} we obtain for each $\alpha >1$ and almost each $x$
it holds $\mathrel{\mathop{\liminf}\limits_{n\rightarrow \infty}}n^{\alpha
}d(T^{n}(x),Y)=\infty .$ Let also suppose that the orbit of $x$ never meet $%
Y $. There is a $c$ such that for all $n$ it holds $n^{\alpha
}d(T^{n}(x),Y)\geq c>0.$

Let us consider the projections $\pi _{q}((q,p))=q,\pi _{p}((q,p))=p.$ Let
us consider an $y$ such that $\forall i\leq n$
\begin{equation}
d(\pi _{q}(T_{\alpha ,\beta }^{i}(x)),\pi _{q}(T_{\alpha ,\beta
}^{i}(y)))\leq \frac{c}{4}i^{-\alpha },~d(\pi _{p}(T_{\alpha ,\beta
}^{i}(x)),\pi _{p}(T_{\alpha ,\beta }^{i}(y)))\leq \frac{c}{4}i^{-\alpha }.
\label{lessa}
\end{equation}%
Then the orbits of $x$ and $y$ are not separated by the discontinuity at the
$n+1$ step. This is true because the orbit of $x$ will stay far away (more
than $\frac{c}{2}i^{-\alpha }$) enough from $Y$ and after the skew
translation $d(\pi _{q}(B(T_{\alpha ,\beta }^{i}(x))),\pi _{q}(B(T_{\alpha
,\beta }^{i}(y))))\leq \frac{c}{2}i^{-\alpha }$ hence when $G_{\alpha }$ is
applied the two points are near enough to avoid to be separated by the
discontinuity.

Now let us estimate the set of points which are near enough to $x$ to
satisfy equation \ref{lessa} after $m$ steps. If
\begin{equation*}
d(\pi _{p}(x),\pi _{p}(y))=d_{p},d(\pi _{q}(x),\pi _{q}(y))=d_{q}
\end{equation*}%
and after $m$ steps, if the orbit of $x$ and $y$ are separated only by the
effect of the skew translation we have that $d_{q}(\pi _{q}(T_{\alpha ,\beta
}^{m}(x)),\pi _{q}(T_{\alpha ,\beta }^{m}(y)))\leq md_{p}+d_{q}$ hence if $%
md_{p}+d_{q}\leq \frac{c}{8}m^{-\alpha }$ the two points are not separated
by the discontinuity at next step. Let us suppose $d_{q}\leq \frac{c}{8}%
m^{-\alpha },$ this gives $d_{p}\leq \frac{c}{4}m^{-\alpha -1}$. By this we
obtain that when $m$ is big enough with respect to $\epsilon $
\begin{equation*}
B(x,m,\epsilon )\supset R=\{y:d(\pi _{q}(x),\pi _{q}(y))\leq \frac{c}{8}%
m^{-\alpha },~d(\pi _{p}(x),\pi _{p}(y))\leq \frac{c}{8}m^{-\alpha -1}\}
\end{equation*}%
the measure of the rectangle $R$ on the right side is $\mu (R)=\frac{c}{8}%
m^{-\alpha }\frac{c}{8}m^{-\alpha -1}=\frac{c^{2}}{64}m^{-2\alpha -1}$ and $%
\alpha $ is near to 1 as wanted. This gives the statement.
\end{proof}

\subsubsection{Logistic map at chaos threshold}

Now we calculate the metric complexity of the orbits of this well
known dynamical system. First let us recall that the Logistic map at
the chaos threshold is a map with zero topological entropy.
Nevertheless the topological complexity of the map $T_{\lambda
_{\infty }}$ is not trivial (see \cite{Ga3}, theorem 22) this means
that the total number of essentially different orbits is not bounded as
time increases. On the contrary as we will see below, the metric
complexity is trivial.

To understand the dynamics of the Logistic map at the chaos threshold let us
use a result of \cite{CE} (Theorem III.3.5.)

\begin{lemma}
\label{teckmann} The logistic map $T_{\lambda _{\infty }}$ at the chaos
threshold has an invariant Cantor set $\Omega $ with the following
properties.

\noindent (1) There is a decreasing chain of closed subsets
\begin{equation*}
J^{0}\supset J^{1}\supset J^{2}\supset \dots ,
\end{equation*}%
each of which contains $1/2$, and each of which is mapped onto itself by $%
T_{\lambda _{\infty }}$.

\noindent (2) Each $J^{i}$ is a disjoint union of $2^{i}$ closed intervals. $%
J^{i+1}$ is constructed by deleting an open subinterval from the middle of
each of the intervals making up $J^{i}$.

\noindent (3) $T_{\lambda _{\infty }}$ maps each of the intervals making up $%
J^{i}$ onto another one; the induced action on the set of intervals is a
cyclic permutation of order $2^{i}$.

\noindent (4) $\Omega =\cap _{i}J^{i}$. $T_{\lambda _{\infty }}$ maps $%
\Omega $ onto itself in a one-to-one fashion. Every orbit in $\Omega $ is
dense in $\Omega $.

\noindent (5) For each $k\in \mathbf{N}$, $T_{\lambda _{\infty }}$ has
exactly one periodic orbit of period $2^{k}$. This periodic orbit is
repelling and does not belong to $J^{k+1}$. Moreover this periodic orbit
belongs to $J^{k}\setminus J^{k+1}$, and each point of the orbit belongs to
one of the intervals of $J^{k}$.

\noindent (6) Every orbit of $T_{\lambda _{\infty }}$ either lands after a
finite number of steps exactly on one of the periodic orbits enumerated in
5, or converges to the Cantor set $\Omega $ in the sense that, for each $k$,
it is eventually contained in $J^{k}$. There are only countably many orbits
of the first type.
\end{lemma}

By this it follows that the metric complexity of this map is trivial, in the
following sense:

\begin{theorem}
In the dynamical system $([0,1],T_{\lambda _{\infty }},\mu )$ if $\mu $ is
some invariant measure supported on the attractor $\Omega ,$ for each $f$ , $%
h_{\mu }^{f}(x)=0$.
\end{theorem}

\begin{proof}
By point 2 of the above lemma \ref{teckmann},  $J^{i}=\cup _{k\leq
2^{i}}J_{k}^{i}$ is the union of $2^i$ intervals, let $\epsilon
_{i}=\max_{k\leq 2^{i}}(diam(J_{k}^{i}))$. By lemma \ref{teckmann},
point 3, if $x,y\in J_{k}^{i}$ then $\sup_{n\geq
0}d(T^{n}(x),T^{n}(y))\leq \epsilon _{i}$. By this we know that for each $%
\epsilon \geq \epsilon _{i}$ and each $x\in J_{m}^{i}$ the set $%
B(x,n,\epsilon )$ contains $J_{m}^{i}$ for each $n.$ Hence $2^{i}$ Bowen
sets are sufficient to cover $J^{i}$ for any $n$. Since the support of the
measure is contained in each $J^{i}$ we have the statement.
\end{proof}

\section{Caracteristic exponents}

The set $B(t,x,\epsilon )$ and its way of shrinking as $t$ increases
describes the initial condition sensitivity of the system around the point $%
x.$

The set will shrink with different speeds at different directions. For
example, the presence of a stable manifold at $x$\ will imply that $%
B(t,x,\epsilon )$ contains for each $n$ a piece of the manifold and
does not shrink in the directions parallel to the manifold. We
introduce a set of numbers $l_{i}$ which describes the shrinking
rate at the different directions. These numbers are in some sense
versions of the positive Lyapunov exponents. In the cases when the
geometry of $B(t,x,\epsilon )$ in nice the numbers $l_{i}$ are
related to the metric complexity, by a result which plays the role
of the Ruelle-Pesin formula.

For simplicity we suppose that $X$ is an open subset of
$\mathbb{R}^{n}$, the case where $X$ is a manifold is similar. Let
us consider the set $S$ of isometries of $\mathbb{R}^{n}.$ Let
$P_{\ell _{1}...\ell _{n}}=[-\frac{\ell
_{1}}{2},\frac{\ell _{1}}{2}]\times ...\times \lbrack -\frac{\ell _{n}}{2},%
\frac{\ell _{n}}{2}]$ be the rectangular parallelepiped with sides
$\ell _{1}...\ell _{n}.$ Let

\begin{equation*}
l_{1}(B(t,x,\epsilon ))=\inf \{\ell _{1}:\exists
~an~isometry~A~s.t.~B(t,x,\epsilon )\subset A(P_{\ell _{1}...\ell _{n}})\}
\end{equation*}

\begin{remark}
\label{remarko} $l_{1}(B(t,x,\epsilon ))$ is a minimum.
\end{remark}

\begin{proof}
This follows by compactness, indeed the space $S$ and the space of all
possible parallelepipeds are locally compact. Moreover, a sequence $%
A_{i}(P_{\ell _{1}^{i}...\ell _{n}^{i}})$ realizing the infimum of $\ell
_{1} $ can be chosen to be a bounded one, hence,by compactness it has a
subsequence having limit.

Since each parallelepiped is compact then this limit parallelepiped will
contain $B(t,x,\epsilon )$, conversely a whole neighborhood of the limit
parallelepiped should not contain $B(t,x,\epsilon )$.
\end{proof}

By this, let us also define%
\begin{equation*}
l_{2}(B(t,x,\epsilon ))=\inf \{\ell _{2}:\exists
~an~isometry~A~s.t.~B(t,x,\epsilon )\subset A(P_{l_{1}\ell _{2}...\ell
_{n}})\}.
\end{equation*}%
By remark \ref{remarko}, $l_{2}$ is well defined, and then more generally we
define $l_{1},...,l_{n}$ as%
\begin{equation*}
l_{i}(B(t,x,\epsilon ))=\inf \{\ell _{i}:\exists
~an~isometry~A~s.t.~B(t,x,\epsilon )\subset A(P_{l_{1},...l_{i-1}\ell
_{i}...\ell _{n}})\}.
\end{equation*}%
Starting from the above defined $l_{1},...,l_{n}$ we can define some
indicator, characterizing the initial condition sensitivity at different
directions.%
\begin{equation*}
\overline{l}_{i}^{f}(x,\epsilon )=\mathrel{\mathop{limsup}\limits_{t%
\rightarrow \infty}}\frac{-log(l_{i}(B(t,x,\epsilon )))}{f(t)},\underline{l}%
_{i}^{f}(x,\epsilon )=\mathrel{\mathop{liminf}\limits_{t\rightarrow \infty}}%
\frac{-log(l_{i}(B(t,x,\epsilon )))}{f(t)}
\end{equation*}%
\begin{equation*}
\overline{l}_{i}^{f}(x)=\mathrel{\mathop{lim}\limits_{\epsilon\rightarrow 0}}%
\overline{l}_{i}^{f}(x,\epsilon ),\underline{l}_{i}^{f}(x)=%
\mathrel{\mathop{lim}\limits_{\epsilon\rightarrow 0}}\underline{l}%
_{i}^{f}(x,\epsilon ).
\end{equation*}%
The numbers $\overline{l}_{i}^{f}(x)$ are in some sense lower
estimations of the way of shrinking of $B(t,x,\epsilon )$ into
different directions. We can also consider the upper estimations
given by

\begin{equation*}
L_{1}(B(t,x,\epsilon ))=\sup \{\ell _{1}:\exists
~an~isometry~A~s.t.~B(t,x,\epsilon )\supset \overset{}{(A(P_{\ell
_{1}...\ell _{n}}))^{\circ }}\}\footnote{%
By $B^{\circ }$ we denote the internal part of $B.$ },
\end{equation*}%
\begin{equation*}
L_{i}(B(t,x,\epsilon ))=\sup \{\ell _{i}:\exists
~an~isometry~A~s.t.~B(t,x,\epsilon )\supset (A(P_{L_{1},...L_{i-1}\ell
_{i}...\ell _{n}}))^{\circ }\},
\end{equation*}%
\begin{equation*}
\overline{L}_{i}^{f}(x,\epsilon )=\mathrel{\mathop{limsup}\limits_{t%
\rightarrow \infty}}\frac{-log(L_{i}(B(t,x,\epsilon )))}{f(t)},\underline{L}%
_{i}^{f}(x,\epsilon )=\mathrel{\mathop{liminf}\limits_{t\rightarrow \infty}}%
\frac{-log(L_{i}(B(t,x,\epsilon )))}{f(t)},
\end{equation*}%
\begin{equation*}
\overline{L}_{i}^{f}(x)=\mathrel{\mathop{lim}\limits_{\epsilon\rightarrow 0}}%
\overline{L}_{i}^{f}(x,\epsilon ),\underline{L}_{i}^{f}(x)=%
\mathrel{\mathop{lim}\limits_{\epsilon\rightarrow 0}}\underline{L}%
_{i}^{f}(x,\epsilon ).
\end{equation*}%
Similar to the traditional Lyapunov exponents the indicators $l_{i}$ and $%
L_{i}$ allows to prove the following inequalities.

\begin{theorem}
If the system is ergodic, it satisfies the assumptions of
proposition \ref{bkg} and the measure $\mu $ is invariant and
absolutely continuous with
bounded density then almost everywhere it holds%
\begin{equation*}
\sum_{i\leq n}\overline{L}_{i}^{f}(x)\geq h_{\mu}^{f}(X,T)\geq \sum_{i\leq n}%
\underline{l}_{i}^{f}(x).
\end{equation*}

\begin{proof}
As before, by proposition \ref{bkg} we have to estimate $\mu (B(t,x,\epsilon
))$ for a typical $x$. We remark that from remark \ref{remarko} it follows
that there is an isometry $A$ such that $B(t,x,\epsilon )\subset
A(P_{l_{1},...l_{n}})$ then $\mu (B(t,x,\epsilon ))\leq \mu
(A(P_{l_{1},...l_{n}}))$.\ Since $\mu $ has bounded density then $\mu
(A(P_{l_{1},...l_{n}}))\leq Const\cdot l_{1}(B(t,x,\epsilon
))l_{2}(B(t,x,\epsilon ))...l_{n}(B(t,x,\epsilon ))$, hence $\log (\mu
(A(P_{l_{1},...l_{n}})))\leq Const_{2}+\log (l_{1}(B(t,x,\epsilon )))+\log
(l_{2}(B(t,x,\epsilon )))+...+\log (l_{n}(B(t,x,\epsilon ))),$ from which,
dividing by $f(t)$ and taking the appropriated limits we obtain \underline{$%
BK$}$^{f}(x)\geq \sum_{i\leq n}\underline{l}_{i}^{f}(x).$ The other
inequality is similar.
\end{proof}
\end{theorem}

\section{Appendix:an example where \protect\underline{$BK$}$^{f}(x,\protect%
\epsilon )\neq \overline{BK}^{f}(x,\protect\epsilon )$\textbf{\label{exa}}}

We will give an example where $\overline{BK}^{f}(x)\neq \underline{BK}%
^{f}(x) $ almost everywhere. For $f(n)=\log (n).$ We remark that by the
Brin-katok theorem such an example is not possible when $f(n)=n.$

Let us consider the two dimensional torus $X=[0,1\ (\hbox{\rm mod\ }%
1)]\times \lbrack 0,1\ (\hbox{\rm mod\ }1)].$ For simplicity, let us equip
it with the $\sup $ distance. If $d^{\prime }$ is the distance on the circle
$[0,1\ (\hbox{\rm mod\ }1)]$ then $d(\left(
\begin{array}{c}
x_{1} \\
y_{1}%
\end{array}%
\right) ,\left(
\begin{array}{c}
x_{2} \\
y_{2}%
\end{array}%
\right) )=\max (d^{\prime }(x_{1},x_{2}),d^{\prime }(y_{1},y_{2})).$ Let us
also define $d_{x}(\left(
\begin{array}{c}
x_{1} \\
y_{1}%
\end{array}%
\right) ,\left(
\begin{array}{c}
x_{2} \\
y_{2}%
\end{array}%
\right) )=d^{\prime }(x_{1},x_{2}),\ d_{y}(\left(
\begin{array}{c}
x_{1} \\
y_{1}%
\end{array}%
\right) ,\left(
\begin{array}{c}
x_{2} \\
y_{2}%
\end{array}%
\right) )=d^{\prime }(y_{1},y_{2}).$

Let us consider $\alpha =0.0505000000000005...=\sum_{n=0}^{\infty }\frac{1}{%
2^{^{2^{2^{n}}}}}$ we have that $\alpha $ is obviously irrational. We define
$T:X\rightarrow X$ as

\begin{equation*}
T=T_{1}\circ T_{2}
\end{equation*}%
where%
\begin{equation*}
T_{1}(x,y)\,=\,(x+\alpha \,\,\,\,\hbox{\rm mod\ }1,\,y)
\end{equation*}%
\begin{equation*}
T_{2}(x,y)\,=\,(x,\,y+\theta (x)\,\,\,\,\hbox{\rm mod\ }1)
\end{equation*}%
where $\theta (q)$ is the discontinuous function over the unit circle
defined in the following way: let us consider the points $\frac{1}{2}$and $%
\frac{1}{2}-\alpha .$ Such points divide the unit circle into two intervals $%
I_{1},I_{2}.$ $\theta (q)=-\frac{1}{4}$ if $q\in I_{1}$ and $\theta (q)=%
\frac{1}{4}$ if $q\in I_{2}$. $T$ at each step rotates on the $x$ direction
and then cuts the torus along the circles $x=\frac{1}{2}$ $x=\frac{1}{2}%
-\alpha $, rotating the torus in opposite directions along the discontinuity
circles. In this system the Lesbegue measure is invariant, hence let us
consider as $(X,T,\mu )$ the system described above with the Lesbegue
measure.

Let us consider the first entrance time of the orbit of $x$ in the ball $%
B(y,r)$ with center $y$ and radius $r$%
\begin{equation*}
\tau _{r}(x,y)=\min (\{n\in \mathbf{N},n>0,T^{n}(x)\in B(y,r)\}).
\end{equation*}%
An irrational $\gamma $ is said to be of type $\nu _{\gamma }$ if
\begin{equation*}
\nu _{\gamma }=sup\{\beta |\mathrel{\mathop{liminf}\limits_{n\rightarrow
\infty}}j^{\beta }(\mathrel{\mathop{min}\limits_{n\in {\bf N}}}|j\gamma
-n|=0)\}.
\end{equation*}%
Lesbegue almost each irrational is of type $1$, but there are
irrationals with type $>1$. For example the $\alpha $ defined above
has type $\infty .$ From the main result of \cite{KS} it can be
deduced that an irrational rotation with angle $\gamma $ of type
$\nu_\gamma >1$ satisfies
\begin{equation}
\underset{r\rightarrow 0}{\lim \sup }\frac{\log \tau _{r}(x,y)}{-\log r}=\nu
_{\gamma }  \label{kimm}
\end{equation}%
for almost each $x$, while
\begin{equation*}
\underset{r\rightarrow 0}{\lim \inf }\frac{\log \tau _{r}(x,y)}{-\log r}\leq
1\ a.e.
\end{equation*}

In other words this implies that for almost each $x$ there are real
sequences $r_{n}$ \ and $r_{n}^{\prime }$ such that $\underset{n\rightarrow
\infty }{\lim }\frac{\log \tau _{r_{n}}(x,\frac{1}{2})}{-\log r_{n}}=\nu
_{\gamma }$ and $\underset{n\rightarrow \infty }{\lim }\frac{\log \tau
_{r_{n}^{\prime }}(x,\frac{1}{2})}{-\log r_{n}^{\prime }}=1.$ Since the
values of $\tau _{r}$ selects times $i$ where the distance $d(T^{i}(x),\frac{%
1}{2})$ is minimal ($d(T^{\tau _{r}(x,y)}(x),\frac{1}{2})=\underset{i\leq
\tau _{r}}{\min }d(T^{i}(x),\frac{1}{2})$ ). This means that there is a
sequence $n_{k}$ such that $d(T^{n_{k}}(x),\frac{1}{2})=\underset{i\leq n_{k}%
}{\min }d(T^{i}(x),\frac{1}{2})\sim n_{k}^{-\frac{1}{\nu _{\gamma }}}$. Now,
coming back to our system we have that $\nu _{\alpha }=\infty $, moreover,
let us remark that $d_{x}(T^{i}(x),\frac{1}{2}-\alpha )=d_{x}(T^{i+1}(x),%
\frac{1}{2})$, hence $\underset{i\leq n_{k}-1}{\min }d_{x}(T^{i}(x),\frac{1}{%
2}-\alpha )\geq \underset{i\leq n_{k}}{\min }d_{x}(T^{i}(x),\frac{1}{2})$.
This means that if the orbit is far from the discontinuity $x=\frac{1}{2},$
then is also far from the other discontinuity. By this let us choose $%
\epsilon <\frac{1}{4}$ and estimate $\mu (B(n,\left(
\begin{array}{c}
x_{0} \\
y_{0}%
\end{array}%
\right) ,\epsilon ))$ where $\left(
\begin{array}{c}
x_{0} \\
y_{0}%
\end{array}%
\right) $ is a typical initial condition satisfying the above equation \ref%
{kimm} with $y=\frac{1}{2}.$ The only source of initial condition
sensitivity is the action of the discontinuities, let us consider the
discontinuity set $Y=\{\left(
\begin{array}{c}
x \\
y%
\end{array}%
\right) \in \mathbb{X}:x=\frac{1}{2}\ or\ x=\frac{1}{2}-\alpha \}$ by
equation \ref{kimm} for each $\delta >0$ there is a sequence $n_{k}$ such
that eventually $n_{k}^{-\delta }=o(\underset{i\leq n_{k}}{\min }%
d(T^{i}(\left(
\begin{array}{c}
x_{0} \\
y_{0}%
\end{array}%
\right) ),Y))$ then for each $\delta >0$ it holds $B(n_{k},\left(
\begin{array}{c}
x_{0} \\
y_{0}%
\end{array}%
\right) ,\epsilon )\supset \lbrack x_{1}-n_{k}^{-\delta
},x_{1}+n_{k}^{-\delta }]\times \lbrack y_{1}-\epsilon ,y_{1}+\epsilon ]$
and $\underset{n\rightarrow \infty }{\lim \inf }\frac{-\log (\mu (B(n,\left(
\begin{array}{c}
x_{0} \\
y_{0}%
\end{array}%
\right) ,\epsilon )))}{\log (n)}=0$ which gives \underline{$BK$}$^{\log
}(\left(
\begin{array}{c}
x_{0} \\
y_{0}%
\end{array}%
\right) )=0.$

For the estimation of $\overline{BK}^{\log }(\left(
\begin{array}{c}
x_{0} \\
y_{0}%
\end{array}%
\right) )$. Let us consider the way the projection on the $x$ circle of the
orbit of the initial point $\left(
\begin{array}{c}
x_{0} \\
y_{0}%
\end{array}%
\right) $ divides the circle. Let us hence consider the sequence $%
x_{0},x_{1}=x_{0}+\alpha ,x_{2}=x_{0}+2\alpha ...$Let us also suppose that
the discontinuity points are not included in the sequence $x_{i}$ (this is
obviously true for a full measure set of initial conditions).

At each time of the form $n_{k}=2^{2^{2^{k}}}$ the unit circle is divided by
the sequence $x_{i}$ into small intervals with length less or equal than $%
\frac{2}{2^{2^{2^{k}}}}.$ This is true because $2^{2^{2^{k}}}$ is the
minimal period of the rotation by the angle $\alpha _{k}=\sum_{n=0}^{k}\frac{%
1}{2^{^{2^{2^{n}}}}}$ and this divides the circle into equal pieces of
length $\frac{1}{2^{2^{2^{k}}}}$. Now
\begin{equation*}
2^{2^{2^{k}}}\sum_{n=k+1}^{\infty }\frac{1}{2^{^{2^{2^{n}}}}}<\frac{1}{%
2^{2^{2^{k}}}}
\end{equation*}%
and then the distance of the first $2^{2^{2^{k}}}$ steps of the two
rotations (with angles $\alpha $ and $\alpha _{k}$) is smaller than the
length of one small interval, giving the required result.

The size of $B(n_{k},\left(
\begin{array}{c}
x_{0} \\
y_{0}%
\end{array}%
\right) ,\epsilon )$ is then estimated by the length of these small
intervals. Indeed we have that, the point $x=\frac{1}{2}$ is contained in
some interval $[x_{i},x_{j}]$. This means that $T^{i}(\left(
\begin{array}{c}
x_{0} \\
y_{0}%
\end{array}%
\right) )$ is on the left of the discontinuity set at a distance less or
equal than $\frac{1}{2^{2^{2^{k}}}},$while $T^{j}(\left(
\begin{array}{c}
x_{0} \\
y_{0}%
\end{array}%
\right) )$ is on the right of the discontinuity line at a distance less or
equal than $\frac{1}{2^{2^{2^{k}}}}.$ This means that
\begin{equation*}
B(n_{k},\left(
\begin{array}{c}
x_{0} \\
y_{0}%
\end{array}%
\right) ,\epsilon )\subset \lbrack x_{1}-n_{k}^{-1},x_{1}+n_{k}^{-1}]\times
\lbrack y_{1}-\epsilon ,y_{1}+\epsilon ],
\end{equation*}
thus $\lim \sup \frac{-\log (\mu (B(n,\left(
\begin{array}{c}
x_{0} \\
y_{0}%
\end{array}%
\right) ,\epsilon )))}{\log (n)}\geq 1$ which gives $\overline{BK}^{\log
}(\left(
\begin{array}{c}
x_{0} \\
y_{0}%
\end{array}%
\right) )\geq 1.$

Since the initial point can be chosen in a full measure set we have

\begin{proposition}
In the above system, for almost each $x$ $\overline{BK}^{\log }(x)\geq 1$
while \underline{$BK$}$^{\log }(x)=0.$
\end{proposition}

\end{document}